\documentclass[english,11pt]{article}
\usepackage{}
\usepackage{amsmath, amssymb, amstext, hyperref}
\usepackage{graphicx, color, framed, soul}
\graphicspath{{./Figures/}}

\title{\bf A Probability Puzzle}

\author{Daniel Otero\footnotemark[1]}

\begin{document}
\maketitle

\newcommand\red[1]{{\color{red}#1}}
\newcommand\blue[1]{{\color{blue}#1}}
\renewcommand{\thefootnote}{\fnsymbol{footnote}}

\footnotetext[1]{School of Engineering and Science, Department of Science, Tecnol\'ogico de Monterrey, Monterrey Campus, M\'exico.}
\footnotetext[2]{Questions, comments, or corrections to this document may be directed to the following email address: {\tt danielotero@tec.mx}.}

\renewcommand{\thefootnote}{\arabic{footnote}}

\begin{abstract}
In this short article, we present a solution to one of the probability puzzles that Daniel Litt, a mathematician at the University of Toronto, posted on his X account earlier this year. The main goal of this note is to show how some of the typical concepts taught in an undergraduate probability course can be used to solve these types of probability problems, which sound simple, but can be very difficult to solve.
\end{abstract}

{\bf Keywords:} Conditional probability, total probability.

\pagestyle{myheadings}
\thispagestyle{plain}
\markboth{D. Otero}{A Probability Puzzle}

\section{Introduction}

In late January, Daniel Litt, an assistant professor at the University of Toronto, posted a probability puzzle on the social network X that sparked lively online interactions among research mathematicians, computer scientists, economists and some math enthusiasts. According to \cite{quanta}, the puzzle was as follows:\\

\textit{``Imagine, Litt wrote, that you have an urn filled with 100 balls, some red and some green. You can't see inside; all you know is that someone determined the number of red balls by picking a number between zero and 100 from a hat. You reach into the urn and pull out a ball. It's red. If you now pull out a second ball, is it more likely to be red or green (or are the two colors equally likely)?''}\\

Of the thousands of people who voted for an answer on X, only about 22\% chose correctly, showing that our probabilistic reasoning can be very off. However, it is important to note that these probability problems can be deceptively difficult and can baffle even great mathematicians \cite{erdos}.\\

In the next section, we present a solution to the puzzle that can be viewed as a conditional probability approach. We would like to emphasize that we use probability concepts learned in an undergraduate probability course, so an understanding of some key ideas should be sufficient to understand the solution. This said, we expect that the proposed approach will be useful for explaining concepts such as conditional probability and total probability at an undergraduate level. In the last section, we conclude with some final remarks.

\section{A Solution}

The proposed solution combines two key concepts: conditional probability and total probability. That is, we will formulate the puzzle as the computation of a conditional probability, and once this expression is established, we will use total probability to compute the terms of the conditional expression. It is worth mentioning that this approach is basically the same as one of the solutions found in \cite{puzzle}, but this solution was obtained independently. In any case, we recommend that you consult \cite{puzzle} if you want to know more about other types of approaches.\\

Let $r_1$ and $r_2$ be the events that you pull out a red ball on the first and second try, respectively. Then, the probability we want to compute is the following:

\begin{equation}
P(r_2|r_1)=\frac{P(r_1\cap r_2)}{P(r_1)}.
\label{prob}
\end{equation}

Now let $R$ be the number of red balls in the urn. As expected, this is our uniform discrete random variable that takes values between 0 and 100. We will use this random variable to compute the probabilities of the latter expression using total probability. Let us begin with the numerator:

\begin{eqnarray*}
P(r_1\cap r_2)&=&\sum_{r=0}^{100}P(r_1\cap r_2|R=r)P(R=r)\\
&=&\sum_{r=2}^{100}P(r_1\cap r_2|R=r)P(R=r).
\end{eqnarray*}
\noindent
Note that the first two terms of the sum are zero since

\[
P(r_1\cap r_2|R=r)=0
\]
\noindent
 when $R=0,1$. As for the denominator, we have that

\begin{eqnarray*}
P(r_1)&=&\sum_{r=0}^{100}P(r_1|R=r)P(R=r)\\
&=&\sum_{r=1}^{100}P(r_1|R=r)P(R=r), ~(P(r_1|R=0)=0).
\end{eqnarray*}
\noindent
Given the latter, Eq. \ref{prob} is equal to

\[
P(r_2|r_1)=\frac{\sum_{r=2}^{100}P(r_1\cap r_2|R=r)P(R=r)}{\sum_{r=1}^{100}P(r_1|R=r)P(R=r)}.
\]
\noindent
Therefore, 

\begin{eqnarray*}
P(r_2|r_1)&=&\frac{\sum_{r=2}^{100}\left(\frac{r}{100}\right)\left(\frac{r-1}{99}\right)\left(\frac{1}{101}\right)}{\sum_{r=1}^{100}\left(\frac{r}{100}\right)\left(\frac{1}{101}\right)}\\
&=&\frac{1}{99}\frac{\sum_{r=2}^{100}(r^2-r)}{\sum_{r=1}^{100}r}\\
&=&\frac{1}{99}\frac{\sum_{r=2}^{100}r^2-\sum_{r=2}^{100}r}{\sum_{r=1}^{100}r}\\
&=&\left(\frac{1}{99}\right)\left(\frac{\frac{(100)(101)(201)}{6}-1-\frac{(100)(101)}{2}+1}{\frac{(100)(101)}{2}}\right)\\
&=&\left(\frac{1}{99}\right)\left(\frac{201}{3}-1\right)\\
&=&\frac{2}{3}.
\end{eqnarray*}
\noindent
Thus, it is more likely to pull out a red ball than a green one given that we picked a red ball in our first draw.\\

By the way, in the previous calculations we used the following formulas, which may look familiar to some readers: the formula of the pyramidal numbers \cite{pyramid}, 

\[
\sum_{i=1}^ni^2=\frac{n(n+1)(2n+1)}{6},
\]
\noindent
and the formula of the triangular numbers \cite{triangle},

\[
\sum_{i=1}^ni=\frac{n(n+1)}{2}.
\]

\section{Final Remarks}

Litt's puzzle is a variation of a classical puzzle known as ``Bertrand's box paradox'' \cite{puzzle}, which is mathematically identical to the ``Monty Hall Problem'' \cite{monty}. A solution to the latter was published by Marilyn vos Savant in Parade magazine in 1990, sparking a heated debate about what Savant had written in her column, nevertheless, Savant was right. Furthermore, this problem has been approached from various perspectives, and it can even be solved using Bayesian networks \cite{ai}, which are a probabilistic model that represents the relationships of a set of random variables by a directed acyclic graph. An implementation of this solution using the Bayesian approach can be found in \cite{dan}.\\

It is possible that Litt's probability puzzle was not as controversial as Savant's solution of the Monty Hall problem, but the collection of problems posed by Litt definitely attracted considerable attention, and even led some probability theorists to do research on topics related to the solutions of some of these brain teasers \cite{doron}.


\newpage

\begin{thebibliography}{1}

\bibitem{quanta}
{\sc E. Klarreich}
{\em Perplexing the Web, One Probability Puzzle at a Time},
in Quanta Magazine (2024).

\bibitem{puzzle}
{\sc Maura B. Paterson and Douglas R. Stinson},
{\em Daniel Litt's Probability Puzzle},
arXiv preprint arXiv:2409.08094 (2024).

\bibitem{erdos}
{\sc Vazsonyi, Andrew},
{\em Which door has the Cadillac},
in Decision Line, 30.1 (1999).

\bibitem{pyramid}
\url{https://en.wikipedia.org/wiki/Square_pyramidal_number}

\bibitem{triangle}
\url{https://en.wikipedia.org/wiki/Triangular_number}

\bibitem{monty}
\url{https://en.wikipedia.org/wiki/Monty_Hall_problem}

\bibitem{ai}
{\sc Russell, Stuart, and Peter Norvig},
{\em AI a modern approach},
Learning 2.3 (2005).

\bibitem{dan}
\url{https://github.com/danotero/MA2014}

\bibitem{doron}
{\sc Shalosh B. Ekhad and Doron Zeilberger},
{\em How to Answer Questions of the Type: If you toss a coin n times, how likely is HH to show up more than HT?},
https://arxiv.org/abs/2405.13561.

\end{thebibliography}
\end{document}